\documentclass[12pt, english]{article}
\usepackage[cp1251]{inputenc}
\usepackage{babel}
\usepackage{amsmath, amssymb,eucal}

\newtheorem{theorem}{\hskip 0.5 cm Theorem}
\newtheorem{lemma}{\hskip 0.5 cm Lemma}

\begin{document}


\title{A STRONG TYPE INEQUALITY FOR CONVOLUTION WITH THE SQUARE ROOT
OF THE POISSON KERNEL}
\author{I.~N.~Katkovskaya, V.~G.~Krotov}
\date{}

\maketitle


\begin{abstract}
The boundary behaviour of convolutions with Poisson kernel and with square root from Poisson kernel
is essentially differs. The first ones have only nontangential limit. For the last ones the convergence is over
domains admittings a logarithmic order of the contact with the boundary (P.Sj\"{o}gren, J.-O.R\"{o}nning).
This result was generalized by authors on the spaces of homogeneous type.

Here we prove the boundedness in $L^p$, $p>1$, and some weighted estimates for the corresponding maximal operator.
Earlier it was known only weak type inequality.

\end{abstract}


\section{Introduction}\label{secIntroduction}


Let
\begin{equation*}
p(z,\theta)=\frac{1}{2\pi}\cdot\frac{1-|z|^2}{|z-e^{i\theta}|^2}
\end{equation*}
be Poisson kernel in the unit circle $B$ of complex plane.

It is well known (see e.g.
\cite{ZygmundTrigonometricSeries1}), that for every function $f\in
L^1[-\pi,\pi]$ the Poisson integral of $f$
\begin{equation*}
Pf(z)=\int_{-\pi}^{\pi}p(z,\theta)f(\theta)\,d\theta
\end{equation*}
converges to $f(\varphi)$ for almost all $\varphi\in[-\pi,\pi]$
provided that $z$ approaches $e^{i\varphi}$ inside a
nontangential domain\footnote{Hereinafter, $a$
means arbitrary fixed positive constant, and we denote by $c$ (with
indexes) we indicate different positive constants depending possibly
on some parameters, which is not important for our purpose.}
\begin{equation}\label{eqNontangentialDomains}
\left\{z:|z-e^{i\varphi}|<a(1-|z|^2)\right\},\;\;a>0.
\end{equation}
J.Littlewood \cite{LittlewoodJLMS27} (see also
\cite{ZygmundTrigonometricSeries1}) has shown that it is the best
result in the following sense. Let $C_0$ be arbitary simple closed
curve through the point $z=1$, which lies totally inside $B$ exept
for this point, and which touches the boundary of $B$ at this
point. Let the curve $C_{\theta}$ be obtained from $C_0$ by rotation
around $z=0$ by the angle $\theta$. Then there is a Blaschke product
which has no limit for almost all $\theta$ when
$z\to e^{i\theta}$ inside $C_{\theta}$.

A study of the convolutions with degrees of the Poisson kernel
\begin{equation}\label{eqPoissonIntegralWithExponent}
P_{l}f(z)=\int_{-\pi}^{\pi}\left[p(z,\theta)\right]^{l+\frac{1}{2}}f(\theta)\,d\theta,\;\;
l\ge0
\end{equation}
was initiated in \cite{Sjogren84}. It is interesting because
$P_{l}(z,\cdot)$ (and $P_{l}f(z)$) satisfy equation
\begin{equation*}
\frac{1}{4}(1-|z|^2)^2\left(\frac{\partial^2u}{\partial
x^2}+\frac{\partial^2u}{\partial y^2}\right)=
\left(l^2-\frac{1}{4}\right)u
\end{equation*}
($\frac{1}{4}(1-|z|^2)^2\Delta u$ is the Laplacian in the hyperbolic
metric). Of course, convergence of $P_{l}f(z)$ can hard be expected without proper normalization.

Let
\begin{equation}\label{eqNormedPoissonIntegralWithExponent}
\mathcal{P}_{l}f(z)=\frac{P_{l}f(z)}{P_{l}1(z)}
\end{equation}
Note that\footnote{Notice $f\asymp g$ means, that there exists a constant
$c>0$, such that $1/c \le f/g\le c$.}
\begin{equation}\label{eqAsymptoticOfPowerOfPoissonKernel}
P_{l}1(z)\asymp \left\{
\begin{array}{cr}
  (1-|z|)^{\frac{1}{2}-l}, & l>0,\\
  \\
  (1-|z|)^{\frac{1}{2}}\log\frac{2}{1-|z|}, & l=0. \\
\end{array}
\right.
\end{equation}
(It doesn't matter what the base of the logarithm is, however it
will be convenient for us to assume the base equals 2).

For $l>0$ boundary behaviour of the integrals
${\mathcal P}_{l}f(z)$ is the same as for
$l=\frac{1}{2}$ which was described above.
The case of $l=0$ is different. For this $l$ the boundary
behaviour of operators (\ref{eqNormedPoissonIntegralWithExponent}) was studied in
\cite{{Sjogren84}}-\cite{Sjogren97}. There it was shown that
${\mathcal P}_{l}f(z)$ converges to $f(e^{i\varphi})$,
for every function $f\in L^p[-\pi,\pi]$ for almost all $\varphi\in[-\pi,\pi]$,
provided that $z$ approaches $e^{i\varphi}$ inside the domain
\begin{equation}\label{eqSjogrenDomains}
\left\{z\in\mathbb{C}:|z-e^{i\varphi}|<
a(1-|z|)\left(\log\frac{2}{1-|z|}\right)^p\right\}.
\end{equation}
The case of $p=1$ is examined in \cite{Sjogren84}, and $p>1$ in
\cite{Roenning93}--\cite{Roenning97}. Note that domains
(\ref{eqSjogrenDomains}) are essentially wider than nontangential
domains (\ref{eqNontangentialDomains}) and admit tangential
approach to the point $e^{i\varphi}$, and the bigger the $p$ there corresponds the higher degree of tangency.

The proof of almost everywhere convergence in
\cite{{Sjogren84}}--\cite{Roenning97} was based on the weak type inequality
\begin{equation}\label{eqWeakTypeInequalityForSquareRoot}
\mu\left\{\mathcal{L}_{p}\left(\mathcal{P}_{0}f\right)>\lambda\right\}\le
c\left(\frac{1}{\lambda}\Vert
f\Vert_{L^p_{\mu}(X)}\right)^p,\;\;\lambda>0,
\end{equation}
for the maximal operator
\begin{equation*}
\mathcal{L}_{p}f(e^{i\varphi})=\sup\left\{|f(z)|: |\arg
z-\varphi|<a(1-|z|)\left(\log\frac{2}{1-|z|}\right)^p\right\},
\end{equation*}
which corresponds to the domains (\ref{eqSjogrenDomains}).

The estimates(\ref{eqWeakTypeInequalityForSquareRoot}) imply almost everywhere tangential convergence for
${\mathcal P}_{l}f(z)$ in a standard way. Furthermore,
in \cite{Roenning93}--\cite{Roenning97} it was shown that the
domains of approach  (\ref{eqSjogrenDomains}) are optimal and
cannot be made wider.

In our paper \cite{KatkovskajaKrotov00} the inequality
(\ref{eqWeakTypeInequalityForSquareRoot}) was extended to the
case of spaces of homogeneous type.

Let $X$ be a compact Hausdorff  space, topology of which is
generated by a quasimetric $d$. This means that the function $d:X
\times X \to [0,\infty )$ satisfies the conditions
\begin{equation}\label{eqQuasiMetricAxiom}
d(x,y)=0\Leftrightarrow x=y,\;\; d(x,y)=d(y,x),\;\;
d(x,y)\leq a_d[d(x,z)+d(z,y)]
\end{equation}
for any $x,y,z\in X$ (the constant  $a_d\ge 1$ does not depend on
the choice of the elements $x,y,z$ in $X$). The family of open balls
\begin{equation*}
B(x,t)=\{y\in X:d(x,y)<t\}
\end{equation*}
forms a base of the topology of $X$. Without loss of generality we can
assume that ${\rm diam}X\le 1$.

Let $\mu$ be a positive Borel measure on $X$ which satisfies
the homogeneity condition
\begin{equation}\label{eqGomoheneuosCondition}
\mu(B(x,t))\asymp t^{\gamma}
\end{equation}
of the order $\gamma>0$ (the constants of the weak equivalence in
(\ref{eqGomoheneuosCondition}) do not depend on $x\in X$ and
$t\in(0, {\rm diam}X]$). The triple $(X,d,\mu)$ is usually called
a space of homogeneous type \cite{CoifmanWeiss}. We denote by
$L^p_{\mu}(X)$, $1\le p<\infty$, the $L^p$-spaces with
respect to $\mu$.

In this paper we study the operators
\begin{equation}\label{eqGeneralP0}
\mathcal{P}_0f(x,t)=\left(\log\dfrac{2}{t}\right)^{-1}
\int_{X}\frac{f(y)}{(d(x,y)+t)^{\gamma}}\,d\mu(y).
\end{equation}

In the particular case of $X=\{z\in\mathbb{C}:|z|=1\}$, $d$ the Euclidean
metric, and $\mu$ the Lebesgue measure (then $\gamma=1$) these
operators essentially coincide with (\ref{eqNormedPoissonIntegralWithExponent})
for $l=0$ (see (\ref{eqAsymptoticOfPowerOfPoissonKernel}) and section \ref{secSomeExamples} below). We are
interested in these operators as limiting case $\alpha=0$ of
the potential type operators in the spaces of homogeneous type
\begin{equation*}
\int_{X}\frac{f(y)}{(d(x,y)+t)^{\gamma-\alpha}}\,d\mu(y),
\end{equation*}
$0<\alpha<\gamma$. The boundary behaviour of such operators have been studied in our paper
\cite{KatkovskajaKrotov00} .

We shall introduce the maximal functions
\begin{equation}\label{eqTangentialMaximalOperators}
\mathcal{L}_{\delta}u(x)=\sup\left\{\vert u(y,r)\vert:\,
d(x,y)<at\left(\log\frac{2}{t}\right)^{\delta}\right\},\;\; a>0,
x\in X,
\end{equation}
depending on the parameter $\delta\ge0$.  We will
use the notation ${\mathcal L}_0=N$ (that's the ``nontangential''
maximal function).

The weak type inequality (\ref{eqWeakTypeInequalityForSquareRoot}) was proved
in \cite{KatkovskajaKrotov00} for the general case of the homogeneous type
spaces
\begin{equation*}
\mu\left\{\mathcal{L}_{\delta}\left(\mathcal{P}_0f\right)>\lambda\right\}\le
c\left(\frac{1}{\lambda}\Vert
f\Vert_{L^p_{\mu}(X)}\right)^p,\;\;\lambda>0,
\end{equation*}
where $\delta=p/\gamma$.
The methods of proof in \cite{KatkovskajaKrotov00} are different from used in \cite{Sjogren84}-\cite{Roenning97}.

From theorem 1 in \cite{KatkovskajaKrotov00} it can deduce also the sharp results
about boundary behaviour of operators (\ref{eqNormedPoissonIntegralWithExponent}) for $l<0$.


\section{The main theorem and its proof}\label{secMainTheorem}


The main goal of this article is the proof that
the inequality (\ref{eqWeakTypeInequalityForSquareRoot})
can be strengthened and replaced by the inequality of a
strong type for $p>1$ even in the general situation.
Namely, the following statement is true.

\begin{theorem}\label{thSrtongTypeInequalityForP0}
If $p>1$ and $\delta=p/\gamma$, then
\begin{equation*}
\|\mathcal{L}_{\delta}\left(\mathcal{P}_0f\right)\|_{L^p_{\mu}(X)}\le
c_p\|f\|_{L^p_{\mu}(X)},
\end{equation*}
where the constant $c_p$  does not depend on $f\in L^p_{\mu}(X)$.
\end{theorem}

Hereinafter, we will put for the simplicity $a=1$ in
(\ref{eqTangentialMaximalOperators}). For the proof we will need
several auxiliary facts. Let us begin with well known ones.

\begin{lemma}\label{lmCovering}
Let $E\subset X$ and $\{B\}$ is any family of balls of bounded
radiuses that cover $E$. Then there is a  finite or enumerable subfamily
$\{B_j\}\subset \{B\}$ such that
\begin{equation}
B_i\cap B_j=\varnothing\, (i\ne j),\;\; E\subset \bigcup _j\rho_d
B_j
\end{equation}
with some constant $\rho_d\ge 1$ that depends only on $d$.
\end{lemma}

Here $\rho B$ is a ball concentric with $B$ whose radius
is $\rho$ times bigger. A proof of the lemma can be found in
\cite{CoifmanWeiss}.

Based on this lemma \ref{lmCovering} and standard techniques
\cite{CoifmanWeiss} one can derive the standard properties
for the Hardy-Littlewood maximal function
\begin{equation*}
Mf(x)=\sup\frac{1}{\mu(B)}\int_{B}|f|\,d\mu,
\end{equation*}
where the supremum is taken over all balls $B=B(y,t)$, containing
the point $x\in X$.

\begin{lemma}\label{lmHardyLittlewoodMaximalFunction}
For any $p\ge 1$ there is a constant $c_p$ such that

1) for all functions $f\in L^1_{\mu}(X)$ and $\lambda>0$
\begin{equation*}
\mu\left\{Mf>\lambda\right\}\le \frac{c_1}{\lambda}\Vert
f\Vert_{L^1_{\mu}(X)},
\end{equation*}

2) for all functions $f\in L^p_{\mu}(X)$
\begin{equation*}
\|Mf\|_{L^p_{\mu}(X)}\le c_p\|f\|_{L^p_{\mu}(X)}.
\end{equation*}
\end{lemma}

Let us consider the parametric family of the approach domains to the boundary
\begin{equation}\label{eqParametricTangentDomain}
 D_{A,\delta}(x)=\left\{(y,t):d(x,y)<t\left(\log\frac{2}{t}\right)^{\delta},\;
A<\left(\log\frac{2}{t}\right)^{\delta}\right\}.
\end{equation}

Note that the second inequality above, determining $D_{A,\delta}(x)$,
is equivalent to following
\begin{equation}\label{eqRestrictionOnT}
 t<\exp\left(1-A^{1/\delta}\right)=\tau_A.
\end{equation}

Let us also define a family of maximal functions
\begin{equation}\label{eqParametricMaximalFunction}
\mathcal{L}_{A,\delta}u(x)=
\sup\left\{\left(\log\frac{2}{t}\right)^{-1}u(y,At): (y,t)\in
D_{A,\delta}(x)\right\},
\end{equation}
which allow to estimate the operator
$\mathcal{L}_{\delta}\left(\mathcal{P}_0f\right)$. This is
 in lemma \ref{lmEsimateP0ByParametricMaximalFunction} below.

In this lemma the following notation
\begin{equation*}
u(y,t)=\frac{1}{\mu(B(y,t))}\int_{B(y,t)}|f|\,d\mu.
\end{equation*}
is used where $f\in L^1_{\mu}(X)$. It is clear that
\begin{equation}\label{eqNontangentialForMeans}
Nu(x)=Mf(x),\;\;x\in X.
\end{equation}

\begin{lemma}\label{lmEsimateP0ByParametricMaximalFunction}
There is a constant $c$ such that
\begin{equation*}
\mathcal{L}_{\delta}\left(\mathcal{P}_0f\right)(x)\le
c\left(Mf(x)+\sum_{\nu=0}^{\infty}\mathcal{L}_{2^{\nu},\,\delta}u(x)\right),
\end{equation*}
for any $x\in X$ and $f\in L^1_{\mu}(X)$.
\end{lemma}

{\bf Proof.} Let $x\in X$ and  $(y,t)\in X\times(0,1)$ satisfies the condition
\begin{equation*}
d(x,y)<\tau=t\left(\log\frac{2}{t}\right)^{\delta}.
\end{equation*}

We split the integral defining the operator (\ref{eqGeneralP0})
into three pieces
\begin{equation*}
\int_{X}\frac{f(z)}{(d(y,z)+t)^{\gamma}}\,d\mu(z)=
\int_{B(y,t)}+\int_{t<d(x,y)\le \tau}+\int_{d(x,y)> \tau}\equiv
I_1+I_2+I_3
\end{equation*}
and we will estimate each of them separately.

First of all we  note that
\begin{equation*}
|I_1|\le u(y,t).
\end{equation*}

Furthermore, if $n=\left[\log_2\frac{\tau}{t}\right]+2$ then
\begin{equation*}
|I_2|\le\sum_{\nu=0}^{n-1}\int\limits_{2^{\nu}t<d(y,z)\le2^{\nu+1}t}\frac{f(z)}{(d(y,z)+t)^{\gamma}}\,d\mu(z)\le
\sum_{\nu=0}^{n-1}\left(2^{\nu}t\right)^{\gamma}\int_{B(y,2^{\nu+1}t)}|f|\,d\mu\le
\end{equation*}
\begin{equation*}
\le c\sum_{\nu=0}^{n-1}\frac{1}{\mu
B(y,2^{\nu+1}t)}\int_{B(y,2^{\nu+1}t)}|f|\,d\mu\le
c\sum_{\nu=1}^{n}u(y,2^{\nu}t).
\end{equation*}

Finally, if $m=\left[\log_2\frac{1}{\tau}\right]+2$ then
\begin{equation*}
|I_3|\le\sum_{\nu=0}^{m-1}\int\limits_{2^{\nu}\tau<d(y,z)\le2^{\nu+1}\tau}\frac{f(z)}{(d(y,z)+t)^{\gamma}}\,d\mu(z)\le
\sum_{\nu=0}^{m-1}\left(2^{\nu}\tau\right)^{\gamma}\int_{B(y,2^{\nu+1}\tau)}|f|\,d\mu\le
\end{equation*}
\begin{equation*}
\le c\sum_{\nu=0}^{m-1}\frac{1}{\mu
B(y,2^{\nu+1}\tau)}\int_{B(y,2^{\nu+1}\tau)}|f|\,d\mu\le
cmMf(x)\le c\log\frac{1}{\tau}Mf(x).
\end{equation*}

Now the conclusion of lemma \ref{lmEsimateP0ByParametricMaximalFunction}
follows automatically from the above inequalities for $I_1$, $I_2$ and $I_3$.

The following lemma plays the key role. The proof follows the ideas from
the paper of the second author \cite{KrotovMIAN89} with suitable
modifications.

\begin{lemma}[main]\label{lmEstimateOfNormParametricMaximalFunction}
For every $p>1$ there is a constant $c_p$ such that
\begin{equation}\label{eqEstimateOfNormParametricMaximalFunction}
\|\mathcal{L}_{A,\delta}u\|_{L^p_{\mu}(X)}\le
c_pA^{-\frac{\gamma}{p}}\|f\|_{L^p_{\mu}(X)}
\end{equation}
for any $f\in L^p_{\mu}(X)$ and $A\ge 1$.
\end{lemma}

{\bf Proof.} Let us consider the Lebesgue sets for the maximal operator
(\ref{eqParametricMaximalFunction})
\begin{equation*}
E_A(\lambda)=\left\{x\in
X:\mathcal{L}_{A,\delta}u(x)>\lambda\right\},\;\;\lambda>0.
\end{equation*}
We shall partition them into parts in the following way. Let $k_{A}$
be the positive integer such that $2^{-k_A-1}<\tau_A\le2^{-k_A}$.
We put
\begin{equation*}
t_{A}(x)=\sup\left\{t<\tau_A:\exists\,(y,t)\in D_{A,\delta}(x),\;
\left(\log\frac{2}{t}\right)^{-1}u(y,At)>\lambda\right\}
\end{equation*}
and for $k\ge k_{A}$ we will define the sets
\begin{equation*}
E_{A,k}(\lambda)=\left\{x\in E_{A}(\lambda):
t_{A}(x)\in\left(2^{-k-1},2^{-k}\right]\right\}.
\end{equation*}
They are measurable and
\begin{equation*}
E_{A,k}(\lambda)\cap E_{A,i}(\lambda)=\varnothing,\;\;
E_{A}(\lambda)=\bigcup_{k=k_{A}}^{\infty}E_{A,k}(\lambda).
\end{equation*}
Also, we will define a modification of the nontangential maximal function
\begin{equation*}
N_Au(x)=\sup\left\{u(y,A\tau):
d(x,y)<\frac{A\tau}{4a_d^2},\;\tau<\tau_A\right\}.
\end{equation*}
It is clear that
\begin{equation}\label{eqInequalityForNAAndN}
N_A(x)\le N(x),\;\;x\in X.
\end{equation}

In the same way as mentioned above let us define the Lebesgue sets
\begin{equation*}
E(\lambda)=\left\{x\in X:N_Au(x)>\lambda\right\},
\end{equation*}
and partition them into
\begin{equation*}
E_{k}(\lambda)=\left\{x\in
E(\lambda):\tau(x)\in\left(2^{-k-1},2^{-k}\right]\right\},
\end{equation*}
where
\begin{equation*}
\tau(x)=\sup\left\{\tau<\tau_A:\exists\,y\;\;d(x,y)<\frac{A\tau}{4a_d^2},\;
u(y,A\tau)>\lambda\right\}.
\end{equation*}
Again, the sets $E_{k}(\lambda)$ are measurable and
\begin{equation*}
E_{k}(\lambda)\cap E_{i}(\lambda)=\varnothing,\;\;
E(\lambda)=\bigcup_{k=k_{A}}^{\infty}E_{k}(\lambda).
\end{equation*}

Let us estimate the measure $\mu E_{A,k}(\lambda)$. For $x\in
E_{A,k}(\lambda)$ there is a pair  $(y_x,t_x)\in
X\times\left(2^{-k-1},2^{-k}\right]$ such that
\begin{equation}\label{eqChooseOfBalls}
d(x,y_x)<t_x\left(\log\frac{2}{t_x}\right)^{\delta},\;\;
u(y_x,At_x)>\lambda\log\frac{2}{t_x}\ge k\lambda.
\end{equation}

Let us consider the family of balls
\begin{equation*}
B_x=B\left(y_x,t_x\left(\log\frac{2}{t_x}\right)^{\delta}\right),\;\;
x\in E_{A,k}(\lambda).
\end{equation*}
By lemma \ref{lmCovering} it is possible to select a finite or countable
subfamily $\{B_{x_j}\}$ with properties
\begin{equation*}
B_{x_j}\cap B_{x_i}=\varnothing,\;\;\mu E_{A,k}(\lambda)\le
c\sum_{j}\mu B_{x_j}.
\end{equation*}

Let us now consider a new family of balls
\begin{equation*}
B^*_{x_j}=B(y_{x_j},At_{x_j}),\;\;j\ge1.
\end{equation*}

We let
\begin{equation*}
\varphi(t)=t\left(\log\frac{2}{t}\right)^{\delta}.
\end{equation*}
It is easy to see that there is a number $k_0\in\mathbb{N}$ such that
\begin{equation}\label{eqInequalityForTangentFunction}
\varphi(t)<\frac{\varphi(\tau)}{4a_d^2}\;\;\text{when}\;\;2^{k_0}t<\tau.
\end{equation}
(sf (\ref{eqQuasiMetricAxiom})).

We will show that the following inclusions are true
\begin{equation}\label{eqMainInclusion}
B^*_{x_j}\subset\bigcup_{i=k-k_0}^{k}E_i(k\lambda),\;\;k\ge
k_0+k_A.
\end{equation}

Let $x\in B^*_{x_j}$ and the pair $(z,\tau)$ be such that
\begin{equation*}
d(x,z)<\frac{A\tau}{2a_d^2},\;\;2^{k_0-k}<\tau<\tau_A.
\end{equation*}
Then by virtue of (\ref{eqQuasiMetricAxiom}), (\ref{eqRestrictionOnT}), (\ref{eqChooseOfBalls}) and
(\ref{eqInequalityForTangentFunction})
\begin{equation*}
d(x_j,z)\le a_d^2\left[d(x_j,y_j)+d(y_j,x)+d(x,z)\right]\le
\end{equation*}
\begin{equation*}
\le
a_d^2\left[t_{x_j}\left(\log\frac{2}{t_{x_j}}\right)^{\delta}+At_{x_j}+\frac{A\tau}{2a_d^2}\right]<
\tau\left(\log\frac{2}{\tau}\right)^{\delta}.
\end{equation*}

Hence, $(z,\tau)\in D_{A,\delta}(x_j)$ (see
(\ref{eqParametricTangentDomain})), however $x_j\in
E_{A,k}(\lambda)$ and $2^{k_0-k}<\tau<\tau_A$, thus
\begin{equation*}
u(z,A\tau)\le\lambda\log\frac{2}{\tau}<\lambda\log2^{k-k_0+1}<k\lambda.
\end{equation*}
This means that $x\notin E_i(k\lambda)$ for $i<k-k_0$.

On the other hand, since $t_{x_j}>2^{-k-1}$, $d(x,y_j)<At_{x_j}$,
the second inequality (\ref{eqChooseOfBalls}) implies that
$x\notin E_{i}(k\lambda)$ for $i>k$. Thus, (\ref{eqMainInclusion}) is proved.

Now using (\ref{eqMainInclusion}) for $k\ge k_0+k_A$ we get
\begin{equation*}
\mu E_{A,k}(\lambda)\le c\sum_{j}\mu B_{x_j}= c\sum_{j}\frac{\mu
B_{x_j}}{\mu B^*_{x_j}}\cdot\mu B^*_{x_j}=
ck^pA^{-\gamma}\sum_{j}\mu B^*_{x_j}=
\end{equation*}
\begin{equation}\label{eqInequalitiesForMeasure1}
=ck^pA^{-\gamma}\mu\left(\bigcup_j B^*_{x_j}\right)\le
ck^pA^{-\gamma}\mu\left(\bigcup_{i=k-k_0}^{k}E_i(k\lambda)\right).
\end{equation}

In a similar but more simple way one can prove the inequlaity
\begin{equation}\label{eqInequalitiesForMeasure2}
\mu\left(\bigcup_{k=k_A}^{k_A+k_0-1}E_{A,k}(\lambda)\right)\le
ck_A^pA^{-\gamma}\mu\left(\bigcup_{i=k_A}^{k_A+k_0-1}E_i(k\lambda)\right).
\end{equation}
Indeed, let
\begin{equation*}
S_A(\lambda)=\bigcup_{k=k_A}^{k_A+k_0-1}E_{A,k}(\lambda)
\end{equation*}
and let $x\in S_A(\lambda)$, then there is a pair $(y_x,t_x)\in
X\times\left(2^{-k_A-k_0},\tau_A\right]$, such that
\begin{equation}\label{eqChooseNewBalls}
d(x,y_x)<t_x\left(\log\frac{2}{t_x}\right)^{\delta},\;\;
u(y_x,At_x)>\lambda\log\frac{2}{t_x}\ge k\lambda.
\end{equation}

In the same way as above we shall consider the set of balls
\begin{equation*}
B_x=B\left(y_x,t_x\left(\log\frac{2}{t_x}\right)^{\delta}\right),\;\;
x\in E_{A,k}(\lambda).
\end{equation*}
According to the lemma \ref{lmCovering} it is possible to
select a finite or countable sub-family $\{B_{x_j}\}$ with
the properties
\begin{equation*}
B_{x_j}\cap B_{x_i}=\varnothing,\;\;\mu S_{A}(\lambda)\le
c\sum_{j}\mu B_{x_j}.
\end{equation*}

Now, let us introduce the new set of balls
\begin{equation*}
B^*_{x_j}=B(y_{x_j},At_{x_j})\;\;(j\ge1),
\end{equation*}
then by virtue of (\ref{eqChooseNewBalls})
\begin{equation*}
B^*_{x_j}\subset\bigcup_{k=k_A}^{k_A+k_0-1}E_{k}(\lambda).
\end{equation*}
By repeating the proof of the inequality (\ref{eqInequalitiesForMeasure1})
we get (\ref{eqInequalitiesForMeasure2}).

Further, using (\ref{eqInequalitiesForMeasure1}) and
(\ref{eqInequalitiesForMeasure2}) one can easily get the following
estimate of the norm
\begin{equation*}
\|\mathcal{L}_{A,\delta}u\|_{L^p_{\mu}(X)}^p=
p\int_0^{\infty}\lambda^{p-1}\mu E_{A}(\lambda)\,d\lambda=
p\sum_{k=k_{A}}^{\infty}\int_0^{\infty}\lambda^{p-1}\mu
E_{A,k}(\lambda)\,d\lambda\le
\end{equation*}
\begin{equation*}
\le c\sum_{k=k_{A}}^{\infty}k^pA^{-\gamma}\sum_{i=k-k_0}^{k}
\int_0^{\infty}\lambda^{p-1}\mu E_i(k\lambda)\,d\lambda.
\end{equation*}
In the last integral we substitute $k\lambda$ instead of
$\lambda$ and get (see also (\ref{eqInequalityForNAAndN}),
(\ref{eqNontangentialForMeans}) and lemma
\ref{lmHardyLittlewoodMaximalFunction})
\begin{equation*}
\|\mathcal{L}_{A,\delta}u\|_{L^p_{\mu}(X)}^p\le
cA^{-\gamma}\sum_{k=k_{A}}^{\infty}
\int_0^{\infty}\lambda^{p-1}\sum_{i=k-k_0}^{k}\mu
E_i(\lambda)\,d\lambda\le cA^{-\gamma}\int_0^{\infty}\lambda^{p-1}
\sum_{i=0}^{\infty}\mu E_i(\lambda)\,d\lambda=
\end{equation*}
\begin{equation*}
=cA^{-\gamma}\int_0^{\infty}\lambda^{p-1} \mu
E(\lambda)\,d\lambda=cA^{-\gamma}\|N_Au\|_{L^p_{\mu}(X)}^p=
cA^{-\gamma}\|Mf\|_{L^p_{\mu}(X)}^p\le
cA^{-\gamma}\|f\|_{L^p_{\mu}(X)}^p.
\end{equation*}
Thus,  lemma \ref{lmEstimateOfNormParametricMaximalFunction} is
proved.

The statement of the theorem \ref{thSrtongTypeInequalityForP0}
follows now directly from a lemmas
\ref{lmHardyLittlewoodMaximalFunction}--\ref{lmEstimateOfNormParametricMaximalFunction}.


\section{Some generalizations and examples}\label{secSomeExamples}


\subsection{Local form of theorem \ref{thSrtongTypeInequalityForP0}.}


First of all we note that the proof of theorem
\ref{thSrtongTypeInequalityForP0} have local character.
This allows us to prove the following generalization of
theorem \ref{thSrtongTypeInequalityForP0}.

\begin{theorem}\label{thLocalForm}
Let $G\subset X$ be an open set and $p>1$. Then for any
compact $K\subset G$ there is a constant $c_p(K)$ such that for
$f\in L^1(X)\cap L^p_{\mu}(G)$ the following inequality is true
\begin{equation*}
\|\mathcal{L}_{\delta}\left(\mathcal{P}_0f\right)\|_{L^p_{\mu}(K)}\le
c_p(K)\left(\|f\|_{L^1_{\mu}(X)}+\|f\|_{L^p_{\mu}(G)}\right).
\end{equation*}
\end{theorem}

{\bf The proof} repeats arguments mentioned above with appropriate
modifications. We describe briefly only relevant changes, which
should be done.

Let $\varepsilon>0$ be small enough number. Let us divide the
operator (\ref{eqGeneralP0}) into two pieces
\begin{equation*}
\mathcal{P}_0f(x,t)=\int_{B(x,\varepsilon)}\frac{f(y)}{(d(x,y)+t)^{\gamma}}\,d\mu(y)+
\int_{X\setminus
B(x,\varepsilon)}\frac{f(y)}{(d(x,y)+t)^{\gamma}}\,d\mu(y).
\end{equation*}
The second integral is estimated as above by
$c\|f\|_{L^1_{\mu}(X)}$ uniformly in $x$ and $t$, and the first one can be
estimated in the same way as it was done above (see lemmas
\ref{lmEsimateP0ByParametricMaximalFunction} and
\ref{lmEstimateOfNormParametricMaximalFunction}). However,
instead of Hardy-Littlewood maximal function $Mf$ one can use "truncated"
maximal function
\begin{equation*}
M_{\varepsilon}f(x)=\sup\frac{1}{\mu(B)}\int_{B}|f|\,d\mu,
\end{equation*}
where the supremum is taken over all balls $B=B(y,t)$ of the
radius $0<t<\varepsilon$, containing the point $x\in X$.


\subsection{Weighted form of theorem \ref{thSrtongTypeInequalityForP0}.}


With the same proof we can to obtain weighted version of
the theorem \ref{thSrtongTypeInequalityForP0}. We will need
some definitions for its statement.

A nonnegative function $\nu $ defined on Borel sets in
$X$ is called an outer measure if it is monotone and subadditive, that is
\begin{equation*}
G_1\subset G_2\;\Rightarrow\; \nu (G_1)\le \nu(G_2),\;\;
\nu \left(\bigcup _{j}G_j\right)\le \sum _j \nu (G_j).
\end{equation*}

If $f$ is a Borel function and $\nu $ is an outer measure on $X$ then we set
\begin{equation*}
\Vert f \Vert_{L_{\nu}^{p}(X)}=\left( p\int ^{\infty }_{0}\lambda^{p-1}
\nu \{ \left| f \right| >\lambda \}\,d\lambda \right) ^{1/p}.
\end{equation*}
For measure $\nu $ this is the usual norm in $L^p_{\nu}(X)$.

\begin{theorem}\label{thWeightedEstimates}
Let $p>1$, $0\le\delta\le p/\gamma$, $\beta=p-\gamma\delta$ and $\nu$ be an
outer measure, satisfying the condition
\begin{equation}\label{eqFrostmanTypeCondition}
\nu(B(x,t))\le ct^{\gamma}\left(\log\frac{2}{t}\right)^{\beta}
\end{equation}
($c$ not depend on $x\in X$ and $t>0$).

Then
\begin{equation*}
\|\mathcal{L}_{\delta}\left(\mathcal{P}_0f\right)\|_{L^p_{\nu}(X)}\le c_p\|f\|_{L^p_{\mu}(X)},
\end{equation*}
where $c_p$ does not depend on $f\in L^p_{\mu}(X)$.
\end{theorem}

{\bf The Proof} word for word copies a proof of the theorem
\ref{thSrtongTypeInequalityForP0}. Only the estimate of the ratio
of measures $\nu B_{x_j}$ and $\mu B^*_{x_j}$ in
(\ref{eqInequalitiesForMeasure1}) and (\ref{eqInequalitiesForMeasure2})
requires to use (\ref{eqFrostmanTypeCondition})  in addition.
Then we obtain
\begin{equation*}
\frac{\nu B_{x_j}}{\mu B^*_{x_j}}\le
\frac{ct^{\gamma}_j\left(\log\frac{2}{t_j}\right)^{\gamma\delta}\left(\log\frac{2}{t_j\left(\log\frac{2}{t_j}\right)^{\delta}}\right)^{\beta}}{(At_j)^{\gamma}}
\le \frac{c\left(\log\frac{2}{t_j}\right)^{\gamma\delta+\beta}}{A^{\gamma}}\le c k^{\gamma\delta+\beta}A^{-\gamma}=
c k^{p}A^{-\gamma}.
\end{equation*}
The rest of the proof passes without change.

Note, that the theorem \ref{thSrtongTypeInequalityForP0} is a particular case
$\beta=0$ of theorem \ref{thWeightedEstimates}.
The last theorem also can give in local form in spirit of the theorem \ref{thLocalForm}.


\subsection{Multidimensional analogues of operator (\ref{eqNormedPoissonIntegralWithExponent}).}


In conclusion, let consider two special cases of operators (\ref{eqGeneralP0}).
 Each of them is a generalization of the operator
(\ref{eqNormedPoissonIntegralWithExponent}) with $l=0$ to a multidimensional case. The theorems
\ref{thSrtongTypeInequalityForP0}-\ref{thWeightedEstimates} can be applied to
both of them.

Let $X=S^{n-1}$ be a unit sphere in $\mathbb{R}^n$, $n\ge 2$,
$\mu$ be the surface  Lebesgue measure on $S^{n-1}$ normalized
by $\mu(S^{n-1})=1$, $d(x,y)=|x-y|$ be the Euclidean metric. Thus,
$\gamma=n-1$ in (\ref{eqGomoheneuosCondition}).

A multidimensional analogue of (\ref{eqPoissonIntegralWithExponent}) be the operator
\begin{equation*}
P_{l}f(x)=\int_{S^{n-1}}\left[p(x,\theta)\right]^{l+\frac{n-1}{n}}f(\theta)\,d\mu(\theta),
\end{equation*}
where
\begin{equation*}
p(x,\theta)=\frac{1-|x|^2}{\left|x-\theta\right|^{n}}
\end{equation*}
is the Poisson kernel for the unit ball (see for example \cite{SteinWeiss}).
In our notation (\ref{eqGeneralP0}) becomes
\begin{equation*}
\mathcal{P}_0f(x)=\left(\log\dfrac{2}{1-|x|}\right)^{-1}
\int\limits_{S^{n-1}}\frac{f(\eta)}{|x-\eta|^{n-1}}\,d\mu(\eta)\asymp
\left(\log\dfrac{2}{t}\right)^{-1}
\int\limits_{S^{n-1}}\frac{f(\eta)}{(|\theta-\eta|+t)^{n-1}}\,d\mu(\eta),
\end{equation*}
where $t=1-|x|$, $\theta=x/|x|$.

Let $X=S^{2n-1}$ be a unit sphere in $\mathbb{C}^n=\mathbb{R}^{2n}$, $\mu$ be
the Lebesgue surface measure, $\mu(S^{n-1})=1$. Let
$d(\zeta,\xi)=\left|1-\langle\zeta,\xi\rangle\right|$ be a nonisotropic
quasimetric ($\langle\cdot,\cdot\rangle$
is the complex scalar product). In this case
$\gamma=n$.

Now it is natural to consider also
the invariant Poisson kernel \cite{RudinFunctionInBall}
\begin{equation*}
P_n(z,\zeta)=\frac{(1-|z|^2)^n}{\left|1-\langle z,\zeta\rangle\right|^{2n}}
\end{equation*}
and by the analogy with (\ref{eqPoissonIntegralWithExponent}) we
arrive at the operators
\begin{equation*}
P_{l}f(z)=\int_{S^{2n-1}}\left[p(z,\eta)\right]^{l+\frac{1}{2}}f(\eta)\,d\mu(\eta).
\end{equation*}
The operator (\ref{eqGeneralP0}) becomes
\begin{equation*}
\mathcal{P}_0f(z)=\left(\log\dfrac{2}{1-|z|}\right)^{-1}
\int\limits_{S^{2n-1}}\frac{f(\eta)}{|1-\langle z,\eta\rangle|^{n}}\,d\mu(\eta)\asymp
\end{equation*}
\begin{equation*}
\asymp\left(\log\dfrac{2}{t}\right)^{-1}
\int\limits_{S^{2n-1}}\frac{f(\eta)}{\left(|1-\langle \zeta,\eta\rangle\right|+t)^{n}}\,d\mu(\eta),
\end{equation*}
with $t=1-|z|$, $\zeta=z/|z|$.

We note in conclusion that other applications of theorems
\ref{thSrtongTypeInequalityForP0}-\ref{thWeightedEstimates} are possible also.
The examples are boundary behaviour of Poisson integrals in polydisk or in Riemannian symmetric
spaces (see \cite{Sjogren88}--\cite{Roenning97}).


\end{document}